\documentclass[a4paper,12pt]{amsart}
\usepackage[labelsep=colon]{caption}
\usepackage{cite}
\usepackage{float}
\restylefloat{figure}
\usepackage{graphicx}
\usepackage{harvard}
\usepackage{hyperref}
\usepackage{mathrsfs}
\usepackage[hang]{subfigure}
\usepackage{tikz}
\usetikzlibrary{arrows,%
  decorations.markings,%
  decorations.pathmorphing,%
  matrix,%
  shapes}

\hypersetup{
  pdftitle={On higher structures},%
  pdfauthor={Nils A. Baas},%
}

\newcommand{\A}{\mathscr{A}}
\newcommand{\C}{\mathscr{C}}
\newcommand{\HH}{\mathscr{H}}
\newcommand{\PP}{\mathcal{P}}

\theoremstyle{definition}

\newtheorem*{definition}{\bfseries Definition}

\theoremstyle{definition}
\newtheorem{example}{\bfseries Example}

\makeatletter
\g@addto@macro\th@definition{\thm@headpunct{:}}
\makeatother

\title{On higher structures}
\author[N.A. Baas]{Nils A.\ Baas$^\ast$}
\dedicatory{Deparment of Mathematical Sciences, NTNU, N-7491
  Trondheim, Norway}
\date{August 31, 2015}
\thanks{$^\ast$Email: \texttt{baas@math.ntnu.no}}

\begin{document}

\begin{abstract}
  In this paper we discuss various philosophical aspects of the
  hyperstructure concept extending networks and higher categories.  By
  this discussion we hope to pave the way for applications and further
  developments of the mathematical theory of hyperstructures.\\

  \noindent \textbf{Keywords:} Higher structures; hyperstructures;
  higher categories; collections of objects.
\end{abstract}

\maketitle

\section{Introduction}
\label{sec:introduction}

In this paper we will discuss the philosophy and ideas we have
developed over the years on higher (order) structures, Baas
\citeyear{EHH,Hnano,HAM,NSCS,NS,SO,HOA}.  Our purpose is to elaborate
on these ideas and introduce new aspects paving the way for further
developments of hyperstructures both in mathematics and other
sciences.

Organization of collections of objects is the clue to create new
structures.  In the hyperstructure concept we have built in the need
to go beyond pairwise interactions and use newly created properties at
one level in the construction of the next level.  See Figure
\ref{fig:formation}.  We will show how this will enable us to
construct new universes and artifical worlds consisting of
hyperstructures.  Our aim is to crystallize the basic ideas in a form
which makes them accessible to workers in many fields.

\section{Organizations}
\label{sec:organizations}

Organizing both abstract and physical matter, organisms, people, etc.\
makes new activities and properties possible.

Hyperstructures --- as we will define them --- pick up and describe
the essence of organizations.  Hence, ``putting a hyperstructure''
onto a situation is an important and useful tool in creating new
activities and structures.  We feel that it is important first to
clarify the main ideas and philosophy behind hyperstructures before
indicating and developing the formal mathematical setting for the
concepts introduced.

Organization matters, as we see for example in:
\begin{itemize}
\item[1)] workers in factories and institutions
\item[2)] energy absorption from light as in photosynthesis --- a
  highly organized process
\item[3)] carbon atoms organized as in diamonds, ashes and graphene
  --- all with different properties
\item[4)] calcium atoms organized in shells and chalk, different types
  of organization with highly different properties.
\end{itemize}

Once we have established a notion of organization in terms of
hyperstructures they may be useful in problem solving, design and data
analysis.

\section{Beyond networks}
\label{sec:beyond_networks}

Networks and graphs play a major role in modelling systems of pairwise
interacting units.  As we have discussed in Baas
\citeyear{NSCS,NS,SO,HOA} there are many examples of systems where
units interact in a collection not reducible to interactions of
subcollections.  For example in Efimov-states we have $3$ particles
interacting, but no pair interactions.  We find the same phenomena in
topological links; see \citeasnoun{NS} and \citeasnoun{BF}.

In a system of interacting units of collections we may ask if there
are further possibilities of interactions when we have exhausted all
pair, triple, $\ldots$ interactions.

Indeed, then we may have collections of units interacting depending on
their properties.  This process may then be extended to collections of
collections, etc.  In order to deal with such situations we need a new
concept --- hyperstructures.

\section{Hyperstructures}
\label{sec:hyperstructures}

We have defined hyperstructures formally in previous papers; Baas
\citeyear{HAM,NSCS,SO}.  Here we will just give an informal definition
which may provide the reader with a better intuition.

Let $\C$ be a collection of objects, and we introduce various bonds or
binding mechanisms of subcollections.  We also assign properties to
subcollections which may play a role in the binding.

\begin{figure}[htbp]
  \centering
  \begin{tikzpicture}[scale=0.5]
    \begin{scope}[xshift=1.25cm]
      \foreach \a in {0,1,...,5}{
        \draw (\a*2,0) circle(0.65cm);
      }

      \node at (12.25,0){$\ldots$};
      \node at (5,-1.5){Basic $0$-order objects};
    \end{scope}

    \begin{scope}[yshift=-3.5cm,xshift=1.25cm]
      \foreach \a in {0,1,...,5}{
        \draw (\a*2,0) circle(0.65cm);
        \draw(\a*2,-0.65) -- (\a*2,-2.3);
        \filldraw[fill=black] (\a*2,-2.3) circle(0.325cm);
      }
      
      \node at (12.25,0){$\ldots$};
      \node at (5,-3.5){Objects with properties};
    \end{scope}

    \begin{scope}[yshift=-9cm,xshift=1.25cm]
      \foreach \a in {0,1,...,5}{
        \draw (\a*2,0) circle(0.65cm);
        \draw(\a*2,-0.65) -- (\a*2,-2.3);
        \filldraw[fill=black] (\a*2,-2.3) circle(0.325cm);
      }

      \draw (0,-2.3) -- (2,-2.3) -- (4,-2.3);
      \draw (6,-2.3) -- (8,-2.3) -- (10,-2.3);
      \node at (12.25,-1.25){$\ldots$};
      \node at (5,-3.5){$1$-bonds};
    \end{scope}

    \begin{scope}[yshift=-15cm]
      \foreach \a in {0,1,2}{
        \draw (\a*2,0) circle(0.65cm);
        \draw(\a*2,-0.65) -- (\a*2,-2.3);
        \filldraw[fill=black] (\a*2,-2.3) circle(0.325cm);
      }

      \draw (-1,-3) rectangle(5,1);
      \draw (0,-2.3) -- (2,-2.3) -- (4,-2.3);

      \draw (1.9,-3) -- (1.9,-4);
      \draw (2.1,-3) -- (2.1,-4);
      \filldraw[fill=black] (2,-4) circle(0.325cm);

      \foreach \a in {4.5,5.5,6.5}{
        \draw (\a*2,0) circle(0.65cm);
        \draw(\a*2,-0.65) -- (\a*2,-2.3);
        \filldraw[fill=black] (\a*2,-2.3) circle(0.325cm);
      }

      \draw (8,-3) rectangle(14,1);
      \draw (9,-2.3) -- (11,-2.3) -- (13,-2.3);

      \draw (10.9,-3) -- (10.9,-4);
      \draw (11.1,-3) -- (11.1,-4);
      \filldraw[fill=black] (11,-4) circle(0.325cm);

      \node at (6.5,-5){$1$-bonds with properties};
      \node at (16,-2.75){$\ldots$};
    \end{scope}

    \begin{scope}[yshift=-22.5cm]
      \foreach \a in {0,1,2}{
        \draw (\a*2,0) circle(0.65cm);
        \draw(\a*2,-0.65) -- (\a*2,-2.3);
        \filldraw[fill=black] (\a*2,-2.3) circle(0.325cm);
      }

      \draw (-1,-3) rectangle(5,1);
      \draw (0,-2.3) -- (2,-2.3) -- (4,-2.3);

      \draw (1.9,-3) -- (1.9,-4);
      \draw (2.1,-3) -- (2.1,-4);
      \filldraw[fill=black] (2,-4) circle(0.325cm);

      \foreach \a in {4.5,5.5,6.5}{
        \draw (\a*2,0) circle(0.65cm);
        \draw(\a*2,-0.65) -- (\a*2,-2.3);
        \filldraw[fill=black] (\a*2,-2.3) circle(0.325cm);
      }

      \draw (8,-3) rectangle(14,1);
      \draw (9,-2.3) -- (11,-2.3) -- (13,-2.3);

      \draw (10.9,-3) -- (10.9,-4);
      \draw (11.1,-3) -- (11.1,-4);
      \filldraw[fill=black] (11,-4) circle(0.325cm);

      \draw (2,-3.9) -- (11,-3.9);
      \draw (2,-4.1) -- (11,-4.1);

      \node at (6.5,-5){$2$-bonds};
      \node at (16,-2.75){$\ldots$};

      \node at (-1,-6){$\vdots$};
      \node at (14,-6){$\vdots$};

      \node at (6.5,-7){etc.};
    \end{scope}
  \end{tikzpicture}
  \caption{Schematic formation of a hyperstructure.  Overlapping may
    occur.}
  \label{fig:formation}
\end{figure}
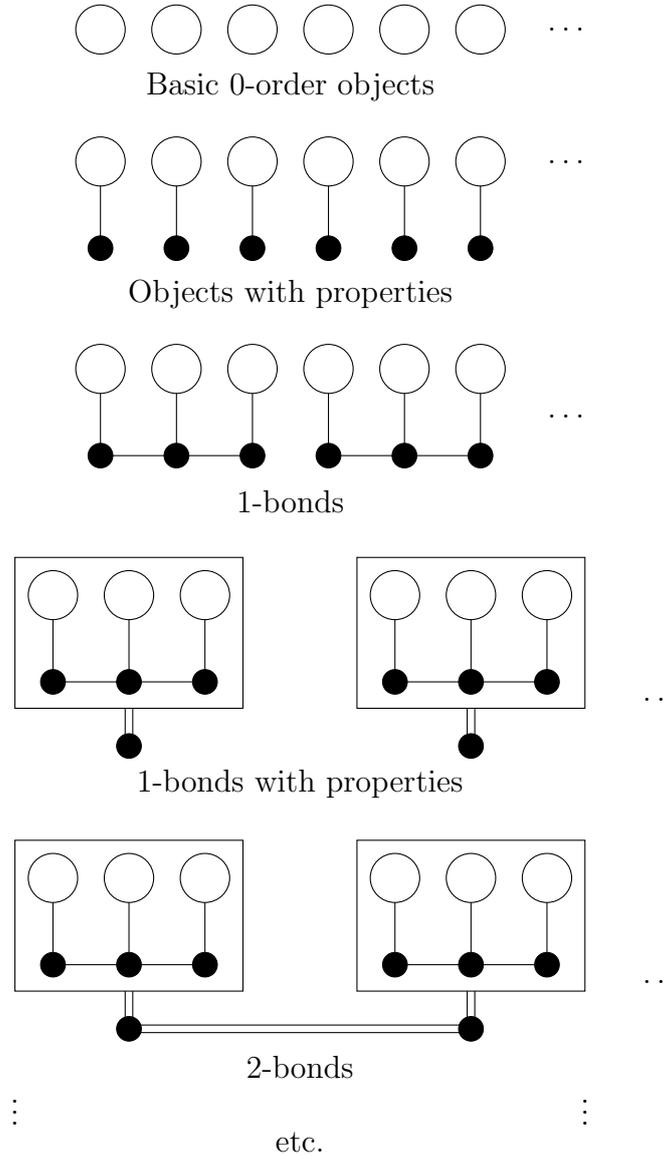

\begin{definition}[informal]
  A \emph{Hyperstructure} $\HH$ consists of a collection of bond sets
  at various levels:
  \begin{equation*}
    \HH = \{B_0,B_1,\ldots,B_n\}
  \end{equation*}
  and the bond levels are connected by maps
  \begin{equation*}
    \partial_i \colon B_{i + 1} \to \PP(B_i)
  \end{equation*}
  taking a bond at level $(i + 1)$ and assigning the collection of
  bonds at level $i$ that it binds.  $B_0$ is the lowest level and
  $B_n$ the top level.  Bonds may be viewed as a kind of general
  multirelations, but it should be noted that properties suitably
  defined at one level play a role in the construction of the next
  level.  Bonds are constructed or created in order to persist as
  functional units according to some presheaf (observer) property.
  See Figure \ref{fig:formation}.
\end{definition}

\begin{figure}[htbp]
  \centering
  \subfigure{
    \centering
    \begin{tikzpicture}[scale=0.7]
      \draw[rotate=15] (0.24396864,-1.4967076) .. controls
      (0.49101374,-2.4657116) and (0.665201,-1.8994577)
      .. (1.5439687,-2.3767076) .. controls (2.4227362,-2.8539574) and
      (3.017664,-3.2517743) .. (4.0039687,-3.4167075) .. controls
      (4.9902735,-3.581641) and (5.629609,-3.3427677)
      .. (6.6239686,-3.2367077) .. controls (7.6183286,-3.1306474) and
      (8.599862,-3.4573336) .. (9.383968,-2.8367076) .. controls
      (10.168076,-2.2160816) and (9.585438,-2.3502257)
      .. (9.943969,-1.4167076) .. controls (10.302499,-0.48318952) and
      (10.588963,-0.26462105) .. (10.743969,0.7232924) .. controls
      (10.898975,1.7112058) and (11.202379,1.6006097)
      .. (10.663969,2.4432924) .. controls (10.125558,3.2859752) and
      (9.734352,3.3049438) .. (8.743969,3.4432924) .. controls
      (7.753585,3.581641) and (5.849801,3.2754364)
      .. (5.0439687,2.6832924) .. controls (4.2381363,2.0911484) and
      (5.064207,1.7455099) .. (4.2239685,1.2032924) .. controls
      (3.3837304,0.661075) and (2.3561227,1.5378368)
      .. (1.5239687,0.9832924) .. controls (0.69181454,0.428748) and
      (-0.003076456,-0.5277036) .. (0.24396864,-1.4967076);

      \begin{scope}[xshift=2cm,yshift=-1.5cm]
        \draw (0,0) .. controls (1,0.5) and (-0.5,1.25) .. (0,2);
        \node[right] at (0.25,0){$Z_1$};
      \end{scope}

      \begin{scope}[xshift=4cm,yshift=-1cm]
        \draw (0,0) .. controls (1,0.5) and (-0.5,1.25) .. (0,2);
        \node[right] at (0.25,0){$Z_2$};
      \end{scope}

      \begin{scope}[xshift=6cm,yshift=-0.5cm,rotate=20]
        \draw[dotted] (0,1) -- (1,1);
      \end{scope}

      \begin{scope}[xshift=8cm,yshift=0cm]
        \draw (0,0) .. controls (1,0.5) and (-0.5,1.25) .. (0,2);
        \node[right] at (0.25,0){$Z_n$};
      \end{scope}

      \node at (7.25,3.5){$Y$};
    \end{tikzpicture}
  } \quad \subfigure{
    \centering
    \begin{tikzpicture}[scale=0.7]
      \draw[rotate=15] (0.24396864,-1.4967076) .. controls
      (0.49101374,-2.4657116) and (0.665201,-1.8994577)
      .. (1.5439687,-2.3767076) .. controls (2.4227362,-2.8539574) and
      (3.017664,-3.2517743) .. (4.0039687,-3.4167075) .. controls
      (4.9902735,-3.581641) and (5.629609,-3.3427677)
      .. (6.6239686,-3.2367077) .. controls (7.6183286,-3.1306474) and
      (8.599862,-3.4573336) .. (9.383968,-2.8367076) .. controls
      (10.168076,-2.2160816) and (9.585438,-2.3502257)
      .. (9.943969,-1.4167076) .. controls (10.302499,-0.48318952) and
      (10.588963,-0.26462105) .. (10.743969,0.7232924) .. controls
      (10.898975,1.7112058) and (11.202379,1.6006097)
      .. (10.663969,2.4432924) .. controls (10.125558,3.2859752) and
      (9.734352,3.3049438) .. (8.743969,3.4432924) .. controls
      (7.753585,3.581641) and (5.849801,3.2754364)
      .. (5.0439687,2.6832924) .. controls (4.2381363,2.0911484) and
      (5.064207,1.7455099) .. (4.2239685,1.2032924) .. controls
      (3.3837304,0.661075) and (2.3561227,1.5378368)
      .. (1.5239687,0.9832924) .. controls (0.69181454,0.428748) and
      (-0.003076456,-0.5277036) .. (0.24396864,-1.4967076);

      \begin{scope}[xshift=2cm,yshift=-1.5cm,scale=0.75,rotate=10]
        \draw (0,0) rectangle(1,3);
        \draw (0,1) rectangle(1,2);
        
        \node[right] at (1.125,0){$Z_1$};
      \end{scope}

      \begin{scope}[xshift=4cm,yshift=-1cm,scale=0.75,rotate=10]
        \draw (0,0) rectangle(1,3);
        \draw (0,1) rectangle(1,2);
        
        \node[right] at (1.125,0){$Z_2$};
      \end{scope}

      \begin{scope}[xshift=6cm,yshift=-0.5cm,rotate=10]
        \draw[dotted] (0,1) -- (1,1);
      \end{scope}

      \begin{scope}[xshift=8cm,yshift=0cm,scale=0.75,rotate=10]
        \draw (0,0) rectangle(1,3);
        \draw (0,1) rectangle(1,2);
        
        \node[right] at (1.125,0){$Z_n$};
      \end{scope}

      \node at (7,3.5){$Y$};
    \end{tikzpicture}
   }
  \caption{$\partial_i Y = Z_i$}
  \label{fig:geometric}
\end{figure}
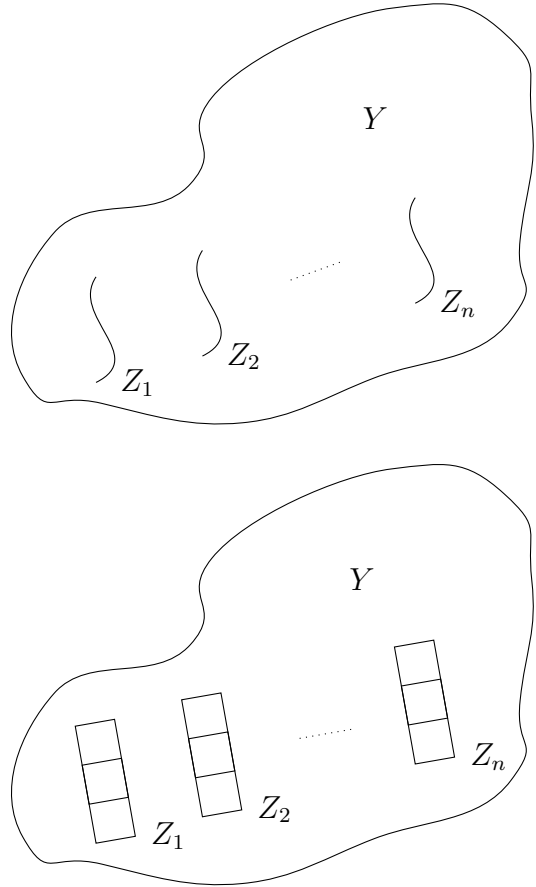

For details see Baas \citeyear{HAM,SO,HOA}, see Figure
\ref{fig:geometric} for geometric illustrations of bonds and the
$\partial_i$'s as ``boundary operators''.  Hyperstructures often
create novel phenomena as discussed in Baas \citeyear{HAM,NSCS,SO} and
may in many contexts be viewed as creators of novelty --- a mechanism
of innovation.  For general structure types a bond of a family of
structures may be a structure containing the family as substructures
of some kind.

We will give some examples of how to build hyperstructures.  Remember
that the hyperstructure concept is an architecture guiding such
constructions where the details can be varied a lot.

\begin{example}[Cones]
  Let the basic set of objects $B_0$ be given as a set of points in
  some Euclidean space.
  \begin{center}
    \begin{tabular}{r@{}c@{\ }p{10cm}}
      $B_0$ & $\colon$ & No property assignment $\Omega_0$.\\
      $B_1$ & $\colon$ & Cones on finite subsets.\\
      $\Omega_1$ & $\colon$ & Cardinality of the coned set.\\[0.2cm]
      $b_1$ & $\colon$ & \begin{tikzpicture}[baseline=-0.95]
        \draw[fill=black] (0,0) -- (120:0.5cm) circle(1pt);
        \draw[fill=black] (0,0) -- (240:0.5cm) circle(1pt);
        \draw[fill=black] (0,0) -- (0:0.5cm) circle(1pt);
        \draw[fill=black] (0,0) circle(1pt);
      \end{tikzpicture} \qquad \begin{tikzpicture}[baseline=-0.95]
        \draw[fill=black] (0,0) -- (40:0.5cm) circle(1pt);
        \draw[fill=black] (0,0) -- (130:0.5cm) circle(1pt);
        \draw[fill=black] (0,0) -- (220:0.5cm) circle(1pt);
        \draw[fill=black] (0,0) -- (310:0.5cm) circle(1pt);
        \draw[fill=black] (0,0) circle(1pt);
      \end{tikzpicture} \quad , \quad $\ldots$\\[0.4cm]
      $\omega_1$ & $\colon$ & $\;\;3$ \qquad\quad $\;\;4$ \qquad ,
                              \quad $\ldots$\\[0.2cm]
      $B_2$ & $\colon$ & Use the property assignment and take for
                         example the cone of cones of just three point
                         sets.  See Figure \ref{fig:cones_cones}.
    \end{tabular}
  \end{center}

  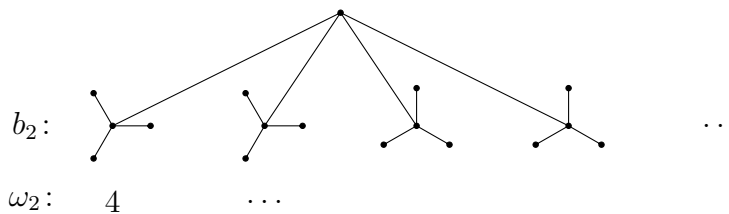
\begin{figure}[htbp]
    \centering
    \begin{tikzpicture}
      \node at (-1,0) {$b_2\colon$};
      \node at (-1,-1) {$\omega_2\colon$};

      \begin{scope}
        \draw[fill=black] (0,0) -- (120:0.5cm) circle(1pt);
        \draw[fill=black] (0,0) -- (240:0.5cm) circle(1pt);
        \draw[fill=black] (0,0) -- (0:0.5cm) circle(1pt);
        \filldraw (0,0) circle(1pt);

        \draw[fill=black] (2,0) --+ (120:0.5cm) circle(1pt);
        \draw[fill=black] (2,0) --+ (240:0.5cm) circle(1pt);
        \draw[fill=black] (2,0) --+ (0:0.5cm) circle(1pt);
        \filldraw (2,0) circle(1pt);

        \draw[fill=black] (4,0) --+ (90:0.5cm) circle(1pt);
        \draw[fill=black] (4,0) --+ (210:0.5cm) circle(1pt);
        \draw[fill=black] (4,0) --+ (330:0.5cm) circle(1pt);
        \filldraw (4,0) circle(1pt);

        \draw[fill=black] (6,0) --+ (90:0.5cm) circle(1pt);
        \draw[fill=black] (6,0) --+ (210:0.5cm) circle(1pt);
        \draw[fill=black] (6,0) --+ (330:0.5cm) circle(1pt);
        \filldraw (6,0) circle(1pt);

        \node at (8,0) {$\ldots$};

        \draw[fill=black] (3,1.5) circle(1pt) -- (0,0);
        \draw (3,1.5) -- (2,0);
        \draw (3,1.5) -- (4,0);
        \draw (3,1.5) -- (6,0);
      \end{scope}

      \begin{scope}
        \node at (0,-1) {$4$};
        \node at (2,-1) {$\ldots$};
      \end{scope}
    \end{tikzpicture}
    \caption{Cones of cones on finite sets}
    \label{fig:cones_cones}
  \end{figure}

  Then we can continue this process by iterated coning but varying the
  number of cones being coned.  For the general picture see Figure
  \ref{fig:cone_cone}.
  \begin{figure}[htbp]
  \centering
  \begin{tikzpicture}
    \draw (0,0) -- (1,1) -- (2,0);
    \draw (0,0) arc (180:360:1cm and 0.15cm);
    \draw[densely dashed] (0,0) arc(180:0:1cm and 0.15cm);

    \begin{scope}[xshift=-0.5cm,yshift=-1.14cm]
      \draw (0,0) -- (1,1) -- (2,0);
      \draw (0,0) arc (180:360:1cm and 0.15cm);
      \draw[densely dashed] (0,0) arc(180:0:1cm and 0.15cm);
    \end{scope}

    \begin{scope}[xshift=-1.5cm,yshift=-3.42cm]
      \draw (0,0) -- (1,1) -- (2,0);
      \draw (0,0) arc (180:360:1cm and 0.15cm);
      \draw[densely dashed] (0,0) arc(180:0:1cm and 0.15cm);
    \end{scope}

    \draw[densely dashed] (0,-1.36) -- (-0.5,-2.32);
  \end{tikzpicture}
  \caption{Several levels of cone structure}
  \label{fig:cone_cone}
\end{figure}
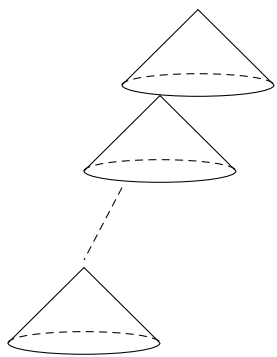
\end{example}

\begin{example}[Manifolds]
  Manifolds are higher dimensional ``surfaces'' looking locally like
  Euclidean spaces.  They are important mathematical objects and occur
  often as parameter and state spaces in applications.

  A possible hyperstructure construction goes as follows.  Let $B_0$
  be a collection of manifolds.
  \begin{center}
    \begin{tabular}{r@{}c@{\ }p{10cm}}
      $\Omega_0$ & $\colon$ & Assign a vector field of some kind for
                              each manifold ($b_0 \mapsto
                              \omega_0$).\\
      $B_1$ & $\colon$ & A bond of a collection of manifolds with vector
                         fields: $b_1\{(b_0^{i_0}, \omega_0^{i_0})\}$
                         is a manifold $b_1$ where the $b_0^{i_0}$'s
                         are embedded as submanifolds and $b_1$ has a
                         vectorfield $\hat{\omega}_1$ extending the
                         $\omega_0^{i_0}$'s.\\
      $\Omega_1$ & $\colon$ & Assigning a vector field $\omega_1$ to
                              each $b_1$.  Possibly $\hat{\omega}_1 =
                              \omega_1$, but this is not required and
                              depends on the situation.\\
      $B_2$ & $\colon$ & Consider a collection of manifolds
                         $\{b_1^{i_1},\omega_1^{i_1}\}$.  A bond $b_2$
                         of this collection will be a manifold with
                         the $b_1^{i_1}$'s embedded and with a vector
                         field $\hat{\omega}_2$ extending the
                         $\omega_1^{i_1}$'s.
    \end{tabular}
  \end{center}
  Then the procedure continues to obtain any number of levels.  See
  Figure \ref{fig:geometric} and \ref{fig:bonds}.
  \begin{figure}[htbp]
    \centering
    \begin{tikzpicture}
      \begin{scope}[rotate=-75,yshift=0.5cm]
        \draw (0,0) circle(2.5cm and 3.5cm);
      \end{scope}

      \draw (-1.5,0.5) to[out=-40,in=75] (-1.25,-1)
        node[below]{$b_0^1$};
      \draw[->] (-1.2,0.1) -- (-0.9,-0.6)
        node[below,xshift=0.35em,yshift=0.45em]{$\omega_0^1$};
      \filldraw (-1.2,0.1) circle(1pt);
      \draw (-0.2,0.8) to[out=-40,in=75] (0.05,-0.7)
        node[below]{$b_0^2$};
      \draw[->] (0.1,0.4) -- (0.4,-0.3)
        node[below,xshift=0.35em,yshift=0.45em]{$\omega_0^2$};
      \filldraw (0.1,0.4) circle(1pt);
      \node[rotate=15] at (1.3,0.7){$\cdots$};
      \draw (2.4,1.4) to[out=-40,in=75] (2.65,-0.1)
        node[below]{$b_0^n$};
      \draw[->] (2.7,1) -- (3,0.3)
        node[below,xshift=0.35em,yshift=0.45em]{$\omega_0^n$};
      \filldraw (2.7,1) circle(1pt);
      \draw[->] (-0.7,1.2) node[left,yshift=-0.15em]{$b_1$} --
        (0.2,1.5) node[right,yshift=0.1em]{$\hat{\omega}_1$};
      \filldraw (-0.7,1.2) circle(1pt);
    \end{tikzpicture}
    \caption{The first level of bonds}
    \label{fig:bonds}
  \end{figure}
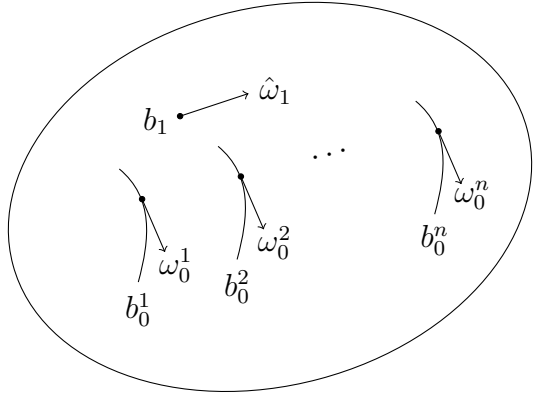
\end{example}

\begin{example}[Organizations]
  Let us consider a collection of agents whose goal is to produce a
  product $P$ with desired property $\omega$.

  One may organize the agents in a hyperstructure as follows.\\
  \begin{center}
    \begin{tabular}{r@{}c@{\ }p{10cm}}
      $B_0$ & $\colon$ & The collection of available agents producing
                         certain products.\\
      $\Omega_0$ & $\colon$ & Assigning a specific product $P_0$ to
                              $b_0\in B_0$.\\
      $B_1$ & $\colon$ & $b_1\in B_1$ binds collections
                         $\{b_0^{i_0},P_0^{i_0}\}$ producing products
                         depending on previous products.\\
      $\Omega_1$ & $\colon$ & assigns a specific product
                              $P_1$ to $b_1$.\\
      $B_2$ & $\colon$ & $b_2\in B_2$ binds collections
                         $\{b_1^{i_1},P_1^{i_1}\}$ producing products
                         depending on previous products.\\
      $\Omega_2$ & $\colon$ & assigns a specific product $P_2$ to
                              $b_2$.\\[0.25cm]
    \end{tabular}
  \end{center}
  This process now continues until the desired $P = P_n$ for some $n$.

  Furthermore, we wanted $P$ to have property $\omega$.  This should
  be tuned in with property assignment at each level (see
  \citeasnoun{HOA}):
  \begin{equation*}
    P_0^{i_0} \mapsto \omega_0^{i_0} \quad , \quad \ldots \quad ,
    \quad P_n \mapsto \omega_n^{i_n}
  \end{equation*}
  and we may then obtain $\omega$ by a globalizer:
  \begin{equation*}
    G \colon \{\omega_0^{i_0}\} \rightsquigarrow \{\omega_1^{i_1}\}
    \rightsquigarrow \cdots \rightsquigarrow \{\omega_n^{i_n}\} \ni
    \omega,
  \end{equation*}
  providing the local to global process through levels.

  In this example a product could mean a physical object with some
  property, but it could also be a solution to a mathematical equation
  with desired properties.  Furthermore, the products could be
  manifolds or stratified spaces with holes (measured by $\Omega$
  being a homology theory) and singularities, properties being vector
  fields of states.
\end{example}

In the examples given here we have simplified the general notation
suppressing the $X_i$'s and $\Gamma_i$'s and taken $B_0$ (the zero
bonds) as the basic objects; see Baas \citeyear{HAM,NSCS}.  For more
examples see also Baas \citeyear{HTD,NSCS,NS,SO}.

In general, if the bonds are given dimensional assignments, higher
bonds will normally have higher dimensions than lower ones.

\section{Local to global}
\label{sec:local_to_global}

In many areas of mathematics there is a need to have methods taking
local information and properties to global ones.  This is mostly done
by gluing techniques using open sets in a topology and associated
presheaves.  The presheaves form sheaves when local pieces fit
together to global ones.  This has been generalized to categorical
settings based on Grothendieck topologies and sites.

The general problem of going from local to global situations is
important also outside of mathematics.  Consider collections of
objects where we may have information or properties of objects or
subcollections, and we want to extract global information.

This is where hyperstructures are very useful.  In \citeasnoun{HOA} we
extended the notion of Grothendieck topologies, sites and
(pre)-sheaves in such a way that ``gluing'' is possible in this
context.

If we are given a collection of objects that we want to investigate,
we put a suitable hyperstructure on it.  Then we may assign ``local''
properties at each level and by the generalized Grothendieck topology
for hyperstructures we can now glue both within levels and across the
levels in order to get global properties.  Such an assignment of
global properties or states we call a \emph{globalizer}.  For a
technical description and details, see \citeasnoun{HOA}.

To illustrate our intuition let us think of a society organized into a
hyperstructure.  Through levelwise democratic elections leaders are
elected and the democratic process will eventually give a ``global''
leader.  In this sense democracy may be thought of as a sociological
(or political) globalizer.  This applies to decision making as well.

In ``frustrated'' spin systems in physics one may possibly think of
the ``frustation'' being resolved by creating new levels and a
suitable globalizer assigning a global state to the system
corresponding to various exotic physical conditions like, for example,
a kind of hyperstructured spin glass or magnet.  Acting on both
classical and quantum fields in physics may be facilitated by putting
a hyperstructure on them.

There are also situations where we are given an object or a collection
of objects with assignments of properties or states.  To achieve a
certain goal we need to change, let us say, the state.  This may be
very difficult and require a lot of resources.  The idea is then to
put a hyperstructure on the object or collection.  By this we create
levels of locality that we can glue together by a generalized
Grothendieck topology.

It may often be much easier and require less resources to change the
state at the lowest level and then use a globalizer to achieve the
desired global change.  Often it may be important to find a minimal
hyperstructure needed to change a global state with minimal resources.

Aspects of this idea were also discussed in \citeasnoun{HOA}.

Again, to support our intuition let us think of the democratic society
example.  To change the global leader directly may be hard, but
starting a ``political'' process at the lower individual levels may
not require heavy resources and may propagate through the democratic
hyperstructure leading to a change of leader.

As discussed in Baas \citeyear{SO,HOA} this is a process which may be
very useful in achieving fusion or fission of sociological, biological
and physical objects.

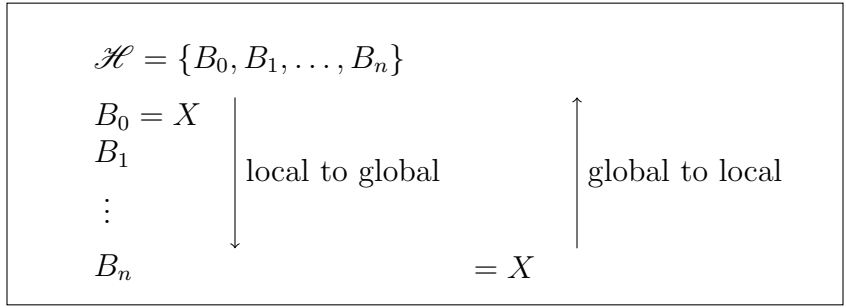
\begin{figure}[htbp]
  \centering
  \begin{tikzpicture}[anchor=west]
    \node at (0,0) {$\HH = \{B_0,B_1,\ldots,B_n\}$};
    \node at (0,-0.75) {$B_0 = X$};
    \node at (0,-1.25) {$B_1$};
    \node at (0,-1.9) {\;$\vdots$};
    \node at (0,-2.75) {$B_n$};

    \draw[->] (2,-0.5) -- node[right]{local to global} (2,-2.5);

    \node at (5,-2.75) {$= X$};

    \draw[->] (6.5,-2.5) -- node[right]{global to local} (6.5,-0.5);

    \draw (-1,0.75) rectangle (10,-3.25);
  \end{tikzpicture}
  \caption{Local to global and global to local processes}
  \label{fig:local_global}
\end{figure}

Hence, hyperstructures facilitates local to global processes, but also
global to local processes.  Often these are called bottom up and top
down processes.  In the global to local or top down process we put a
hyperstructure on an object or system in such a way that it is
represented by a top level bond in the hyperstructure.  This means
that to an object or system $X$ we assign a hyperstructure
\begin{equation*}
  \HH = \{B_0,B_1,\ldots,B_n\}
\end{equation*}
in such a way that $X = b_n$ for some $b_n \in B_n$ binding a family
$\{b_{n - 1}^{i_1}\}$ of $B_{n - 1}$ bonds, each $b_{n - 1}^{i_1}$
binding a family $\{b_{n - 2}^{i_2}\}$ of $B_{n - 2}$ bonds, etc.\
down to $B_0$ bonds in $\HH$.  Similarly for a local to global process.
To a system, set or collection of objects $X$, we assign a
hyperstructure $\HH$ such that $X = B_0$, (Figure
\ref{fig:local_global}).

A hyperstructure on a set (space) will create ``global'' objects,
properties and states like what we see in organized societies,
organizations, organisms, etc.  The hyperstructure is the ``glue'' or
the ``law'' of the objects. In a way, the globalizer creates a kind of
higher order ``condensate''.

Hyperstructures represent a conceptual tool for translating
organizational ideas like for example democracy, political parties,
etc.\ into a mathematical framework where new types of arguments may
be carried through.  This metaphor may be useful in other scientific
contexts as well, like engineering, architecture, etc.

In purely mathematical situations one often wants to form global
solutions of a problem (like differential equations) from local ones.
In such situtations hyperstructures may be useful in introducing bonds
and levels.  Sometimes they may ``force'' a solution to a problem or
it may come ``naturally'' in a compatible way.

\section{The general binding problem}
\label{sec:the_general_binding_problem}

An important special case of local to global situations that we have
discussed, is the so-called binding problem in the neurosciences.  It
is often being illustrated by the ``grandmother problem'': how do you
remember your grandmother?  Such a memory may require many types of
information: visual, auditory, olfactory, etc.\ from different parts
of the brain.  How is it all being integrated into one \emph{global}
meaning?  

In his study D.\ Hebb \citeyear{Hebb} introduced cell assemblies and
hierarchies of assemblies in order to study this problem; see also A.\
Scott \citeyear{Scott}.  Basically in our terminology Hebb was putting
a hyperstructure on the collection of neurons (in the brain) that he
studied.  As we have pointed out, the essence of interesting
organizations is given by a hyperstructure.  See also
\citeasnoun{Cog}.

Therefore our local to global problem is a \emph{generalized binding
  problem} formulated in a precise way.  The existence of a
\emph{globalizer} is a solution of the given binding problem.

We will hope that this could be a useful framework for studying the
binding problem in the brain in the future, and that the analysis of
experimental data may lead to a quantitive hyperstructure.

\section{Higher entanglement}
\label{sec:higher_entanglement}

\begin{figure}[htbp]
  \centering
  \includegraphics[scale=0.5]{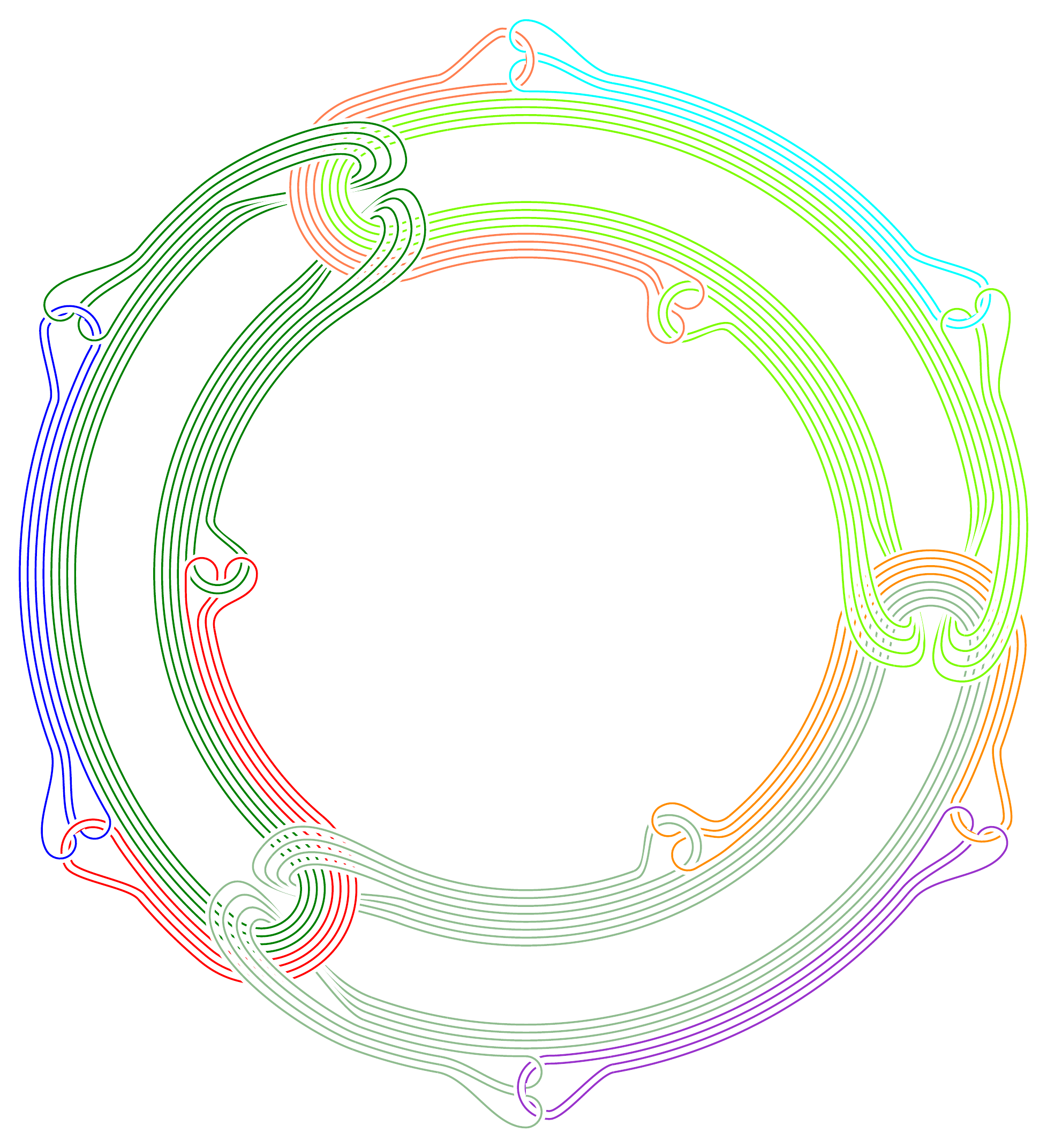}
  \caption{A second order Brunnian ring.}
  \label{fig:brunnian}
\end{figure}

In topology we study entanglement of links in three dimensional space
and of manifolds in higher dimensions.  In \citeasnoun{NS} we studied
links and their higher order versions.  Brunnian rings represent
interesting examples.  We can use their ring property to iterate the
construction process and form second order Brunnian rings etc.  See
Figure \ref{fig:brunnian}.

This gives higher forms of entanglement and forms interesting
hyperstructures of collections of rings or circles.  The entanglement
of Brunnian rings (``remove one ring and they all fall apart'')
corresponds by an interesting analogy to quantum mechanical states of
ultracold gases where $3$ (or $4$) particles may interact but with no
pair interactions.

In \citeasnoun{NS} we therefore conjectured that any form of higher
topological entanglement should have an analogue of quantum mechanical
states.

Furthermore, this is part of a much more general scheme: representing
particles by lowest level bonds (or objects) in a hyperstructure.
Then the higher level bonds will suggest new types of interactions of
the particles --- including higher interactions of higher collections
of particles; see Baas \citeyear{NS,HOA} and \citeasnoun{BF}.

An interesting question is then: how many of these are realizable in
the physical universe?  Instead of starting with various types of
links one may start with graphs, polyhedra, solids, manifolds, ravels,
etc.\ (see \citeasnoun{NS}, \citeasnoun{CEH} and \citeasnoun{C}) and
form various forms of entanglement.  Then, inspect their properties
and use these enforcing new second order tangles, etc.\ forming a
hyperstructure.

All these higher geometric forms may act as models for designing and
synthesizing new types of molecules and materials, as suggested in
\citeasnoun{NSCS}, \citeasnoun{BS} and \citeasnoun{BSS}.  These new
``higher materials'' will probably have new emergent properties --- a
subject for further study.  As pointed out in \citeasnoun{NSCS}
hyperstructured materials may be interesting in the search for high
temperature superconducters.  Hyperstructures offer a variety of
designs for new types of ``higher'' systems and materials.

Hyperstructures may also act as scaffolds in both theoretical and
physical/chemical processes.  An example illustrates this:

Form a Brunnian ring of some length and order by DNA-molecules.  This
is a well-known procedure.  Attach peptides to the DNA rings.  When
they have formed a topological copy, cut the DNA rings, and they fall
apart (the Brunnian property) and we are left with a Brunnian polymer.

Our description of how to induce higher interactions on a collection
of particles from a given hyperstructure is also a kind of scaffolding
process.  Putting a hyperstructure on a situation may sometimes be
thought of as a useful scaffold.

\section{New universes}
\label{sec:new_universes}

Hyperstructures capture the essence of organized structures, and are
very rich with respect to compositional laws (interactions) ---
through ``gluing'' at several levels.

Hence, we may form rich universes where the constituent objects or
``organisms'' are hyperstructures of some kind --- a society of
hyperstructures.  Out of these we may build artifical worlds like in
A-life; see Baas \citeyear{HAM,SO,HOA}.  The point is that this gives
a lot of freedom to do both theoretical and computer experiments.  We
call these $\HH$-worlds or $\HH$-universes, letting now the prefix
$\HH$ stand for hyperstructure(d).  They will often be hyperstructures
of hyperstructures$\ldots$

\emph{Abstract matter} (Baas \citeyear{HAM,SO,HOA}) is given by the
collections formed in these universes --- both statically and
dynamically --- as hyperstructures by the given principles.  Here the
laws and mechanisms can be freely varied to form structures with
desired properties or just to try out new mechanisms and structures.

These universes are really a new kind of testing reactors for
realizing ``fantasy'' structures.  Forces and interactions are
represented by bonds and bond structures.  In the abstract matter
universes one can make new structures that never existed in other
universes --- abstract, mental or physical.  This again may inspire
constructions in our real physical universe.

For example, we may view an organism as an $\HH$-world of
$\HH$-structured cells.  We may for example assume that
$\HH$-organisms are regulated by an $\HH$-genome inducing
$\HH$-structures on the ``brain'', ``immune system'', etc.  This leads
to interesting questions regarding brains: What kind of
$\HH$-structures are there?  What kind of entanglements or Brunnian
phenomena?

As suggested in Baas \citeyear{HAM,SO,HOA}, failures in the
$\HH$-genome structure may lead to ``artifical or abstract cancer''
which then can be studied in these abstract organisms.

Also on the mathematical side we may introduce interesting types of
$\HH$-structures to be studied:
\begin{center}
  \begin{tabular}{r@{}c@{\ }p{9cm}}
    $\HH$-shapes & $\colon$ & shapes of shapes, $\ldots$\\
    $\HH$-forms & $\colon$ & forms of forms, $\ldots$\\
    $\HH$-languages & $\colon$ & languages of languages,$\ldots$,
                                 forming a ``brain'' for higher
                                 dimensional geometric languages,
                                 where given an alphabet $\A$,
                                 $\HH(\A)$ is what the language can
                                 express\\
    $\HH$-fields & $\colon$ & fields of fields,$\ldots$, quantum
                              fields of quantum fields, $\ldots$\\
    $\HH$-materials & $\colon$ & materials of materials, $\ldots$\\
    $\HH$-geometries & $\colon$ & geometries of geometries,$\ldots$,
                                  symmetries of symmetries, $\ldots$\\
    $\HH$-equivalences & $\colon$ & equivalence classes of equivalence
                                    classes, $\ldots$\\
    $n\HH$-structures & $\colon$ & hyperstructures of hyperstructures
                                   of $\ldots$ ($n$ times)\\
    $\HH$-engineering & $\colon$ & $\HH$ structured engineering
                                   products and designs
  \end{tabular}
\end{center}
\noindent are all interesting notions.

$\HH$-geometries may be introduced in a similar way to
$\HH$-Grothendieck topologies in the sense that we may have Riemannian
structures at each level in a compatible way.  We may then make
geometries where we change the laws of nature similar to particle
attraction in energy wells.\sloppy

In such universes we can play all sorts of games with hyperstructures
and vary the rules as we wish.  We may also choose whether we want to
start out with simple or complex objects as a basis for building the
hyperstructures.

In static universes of hyperstructures we may also introduce dynamics
via compositions of bonds or other rules --- again hyperstructured
rules.  This may change or produce new hyperstructures for example by
fusion, fission or other similar $\HH$-processes.  In an environment
(of $\HH$-type) natural selection may occur and evolution into new
types of hyperstructures may take place.

The moral is that ``everything'' being organized in the universe is a
hyperstructure of some kind.  We just have to identify them.  In any
universe hyperstructures are basic tools in organizing constructions
and in the study of collections of objects.

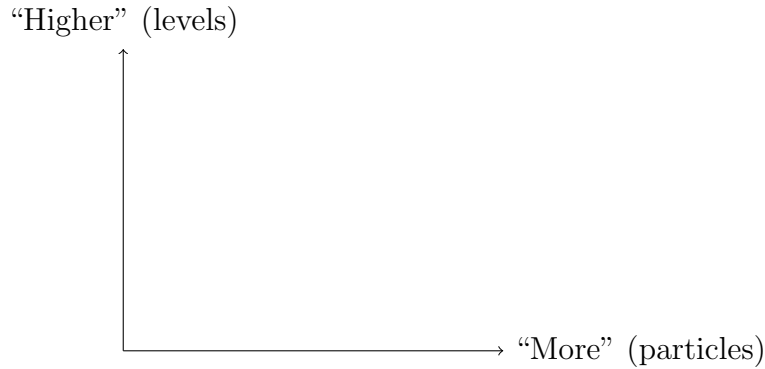
\begin{figure}[htbp]
  \begin{center}
    \begin{tikzpicture}
      \draw[->] (0,0) -- (5,0) node[right]{``More'' (particles)};
      \draw[->] (0,0) -- (0,4) node[above]{``Higher'' (levels)};
   \end{tikzpicture}
  \end{center}
  \caption{``$\HH$-organization is different''.}
  \label{fig:height}
\end{figure}

In \citeasnoun{Anderson} it is argued that ``more is different'',
meaning that in a system of particles new properties emerge as the
number increases.  We would like to add that ``higher is different'',
meaning that new properties emerge with the height of the order in the
hyperstructures; see Figure \ref{fig:height}.

Hyperstructures offer a general organizing principle and are the basic
tools in constructing new types of universes --- abstract or physical.

\section{Evolution}
\label{sec:evolution}

Much of the intuition around hyperstructures comes from thinking of
them as evolutionary structures.  They are designed and defined in the
same way as evolution works: collections interact forming new bonds of
collections with new properties, these being selected for further
interactions forming the next level of bonds etc.  In a sense nature or
the environment acts as a kind of observer (or ``observation
sheaf'').  The success of evolutionary structures makes their
theoretical counterparts --- hyperstructures --- a useful design
model.

In biology we notice that organisms (real hyperstructures) are very
efficient energy devices both with respect to absorption and
dissipation.  Hence, it may turn out that hyperstructures are useful
designs for many types of energy devices.

A fundamental goal is to design and construct devices producing energy
such that the ratio
\begin{equation*}
  \frac{\text{Energy -- Out}}{\text{Energy -- In}}
\end{equation*}
is as large as possible.  Inspired by biological systems it seems
natural to suggest that suitably hyperstructured devices might do
this.  See also \citeasnoun{P} in this direction.  Nuclear fusion may
also be facilitated through a hyperstructured process.

Often when we start out with some basic interacting structures or
objects, these will in many cases have an inbuilt drive to increase
the complexity in order to handle certain situations, solve problems,
etc.  This leads to an evolutionary process and building up
hyperstructures describing higher structures in general.  In such a
process of ``evolving'' structures of higher order one may fuse some
by new bonds and let others dissolve or undergo fission by taking
``boundaries''.

Evolution shows that organization matters and that hyperstructures are
often favoured.  Fusion of for example genes gives new opportunities
and new properties.

\section{Hyperstructured stacking of shapes}
\label{sec:stacking}

Often higher orderness and high dimensions are connected, but in the
synthesis of physical and chemical structures we are limited by our
three dimensional space.  We will elaborate on a procedure given in
\citeasnoun{SO} on how to produce higher geometric and physical
structures by packing higher dimensions into three space.

Suppose that we are given an alphabet of objects in three space
(circles, rings, surfaces, molecules,$\ldots$) and that we have
mechanisms for connecting (binding) them; see Figure
\ref{fig:connecting}.

\begin{figure}[htbp]
  \centering
  \subfigure[1D and 2D objects.]{
    \centering
    \begin{tikzpicture}[scale=0.8,font=\scriptsize]
      \foreach \l in {0,0.5,1,1.5}{
        \draw (0,\l) -- (2,\l);
      }

      \node at (3.75,0.75){(bound together)};
      \node at (7.5,0.75){may be represented by};

      \draw (10,0) rectangle(12,1.5);

      \node at (13.5,0.75){(a rectangle)};
    \end{tikzpicture}
  } \quad \subfigure[New chains being formed by 2D and 3D objects.]{
    \centering
    \begin{tikzpicture}[scale=0.8]
      \draw (0,0) rectangle (2,2);

      \draw (0.5,2.5) -- (2.5,2.5) -- (2.5,0.5);
      \draw (0.5,2.5) -- (0.5,2.0675);
      \draw[dashed] (0.5,1.875) -- (0.5,0.5) -- (1.875,0.5);
      \draw (2.0675,0.5) -- (2.5,0.5);

      \draw (1,3) -- (3,3) -- (3,1);
      \draw (1,3) -- (1,2.5675);
      \draw[dashed] (1,2.375) -- (1,1) -- (2.375,1);
      \draw (2.5675,1) -- (3,1);
      \draw[dotted] (1,3) -- (0.5,2.5);
      \draw[dotted] (3,3) -- (2.5,2.5);
      \draw[dotted] (3,1) -- (2.5,0.5);

      \path (5,1.25) edge[->,bend right=-30]
        node[midway,above,yshift=0.25cm]{may be represented by $3$
          dim box} (9,1.25); 

      \draw (11,0) rectangle (13,2);
      \draw (11,2) -- (12,3) -- (14,3) -- (14,1) -- (13,0);
      \draw (13,2) -- (14,3);
      \draw[dashed] (11,0) -- (12,1) -- (14,1);
      \draw[dashed] (12,3) -- (12,1);
    \end{tikzpicture}
  } \quad \subfigure[Binding structures.]{
    \centering
    \begin{tikzpicture}[scale=0.8]
      \begin{scope}
        \foreach \y in {0,0.75,1,,1.25,1.5}{
          \draw (0,\y) node[left,yshift=-0.05cm]{$\cdots$} -- (1,\y);
        }

        \node at (-0.375,0.45){$\vdots$};
        \node at (0.5,0.45){$\vdots$};
        \node[below,yshift=-0.15cm] at (0.5,0){$b_0$};
      \end{scope}

      \begin{scope}[xshift=4cm]
        \foreach \y in {0,0.75,1,,1.25,1.5}{
          \draw (0,\y) -- (1,\y);
        }
        
        \draw (0,0) -- (0,1.5);
        \node at (0.5,0.45){$\vdots$};
        \node[below,yshift=-0.15cm] at (0.5,0){$b_1$};
      \end{scope}

      \begin{scope}[xshift=8cm]
        \foreach \y in {0,0.75,1.125,1.5}{
          \draw (0,\y) -- (1,\y);
          \foreach \x in {0.1,0.35,0.6,0.85}{
            \draw (\x,\y) -- (\x+0.15,\y+0.15);
          }
        }

        \draw (0,0) -- (0,1.5);
        \node at (0.5,0.5){$\vdots$};
        \node[below,yshift=-0.15cm] at (0.5,0){$b_2$};
      \end{scope}

      \begin{scope}[xshift=12cm]
        \foreach \y in {0,0.75,1,1.25,1.5}{
          \draw[decorate,decoration=zigzag] (0,\y) -- (1,\y);
        }

        \draw (0,0) -- (0,1.5);

        \node at (0.5,0.5){$\vdots$};
        \node[below,yshift=-0.15cm] at (0.5,0){$b_n$};
      \end{scope}

      \foreach \x in {2,6,10}{
        \path (\x,1) edge[->,bend right=-30] (\x+1,1);
      }
    \end{tikzpicture}
    \centering
  }
  \caption{}
  \label{fig:connecting}
\end{figure}

Then we connect them into a line or curve, a bar or a solid cylinder.
This is similar to making polymers, basically $1$D objects.  Next we
take these bars and connect them into a flat surface (rectangle), and
proceed in the same way to form cubes (solids) --- possibly with holes
and other topological and geometric properties.  Then we have
exhausted our three dimensions --- what next?

We start over again and form lines (curves) of the newly constructed
``cubes'', and through this way of ``stacking'' geometric structures
we may synthesize forms (topologies and geometries) of higher
dimensions.  The stacking is determined by the type of bonds
available.  We call this procedure \emph{Hyperstructured Stacking of
  Shapes}.  For interesting aspects of stacking of graphene, see
\citeasnoun{Gibney} and references therein.

Stacking is a counterpart to clustering of objects (including
hierarchical clustering) where we just form zero dimensional geometric
objects.

Persistent homology may be useful in detecting and identifying higher
structures.

This procedure is similar to the way we describe higher cobordism
categories:
\begin{itemize}
\item[] objects as points --- $0$-morphisms ($0$-bonds)
\item[] points generate intervals (curves) --- $1$-morphisms ($1$-bonds)
\item[] intervals generate surfaces --- $2$-morphisms ($2$-bonds)
\item[] \hspace*{1cm}$\vdots$
\item[] $(n - 1$)-manifolds generate $n$-manifolds --- $n$-morphisms
  ($n$-bonds).
\end{itemize}
See for example \citeasnoun{L}.

In general $n$-dimensional objects like manifolds, stratified spaces,
polyhedra, etc.\ may be formed like hyperstructures starting by points
($0$ dimensional objects), using these to generate $1$ dimensional
objects, and using these again to generate $2$ dimensional objects,
etc.

\section{Beyond categories}
\label{sec:beyond_categories}

In mathematics most work on higher structures has been in higher order
logic, higher set theory and higher categories and their applications.
In recent years the theory of higher categories has undergone a strong
development, and is certainly a model for what we would also like to
do with hyperstructures.  As we have pointed out in several papers,
Baas \citeyear{NSCS,NS,SO,HOA}, categories are based on morphisms
which are in our terminology bonds of pairs --- source and target.

This is the picture we extend in hyperstructures replacing morphisms
by bonds of many objects --- like a manifold (cobordism) connecting
(binding) several boundary components.  In addition we assign
properties to bonds and these are being used in forming the next level
of bonds.  Furthermore, once a hyperstructure has been constructed we
often want to study how local properties and states may be glued
together to a global property or state.

In \citeasnoun{HOA} we discussed how this can be done by generalized
Grothendieck topologies, presheaves and globalizers.  Many interesting
questions regarding the mathematical theory of hyperstructures remain
to be explored.  Categories and their higher versions are all based on
``pair interactions''.  This is also the case with simplicial sets and
simplicial structures since everything is basically generated by edges
defined by pairs of vertices.  Abstractly we define a simplex of a
simplicial complex to have the property that any (non-empty) subset is
a simplex as well.  We think of simplices as bonds of its faces.
Within the bond picture it would be natural to weaken this to an
$\HH$-(hyperstructured) complex where an $\HH$-simplex would come as a
family of subsets being the sub-$\HH$-simplices.

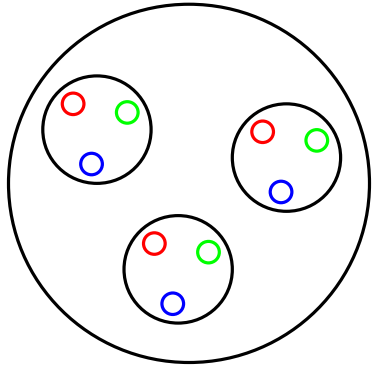
\begin{figure}[htbp]
  \centering
  \begin{tikzpicture}[scale=0.95]
    \draw[very thick] (0,0) circle(2.5cm);

    \begin{scope}[scale=0.3,xshift=-4.25cm,yshift=2.5cm]
      \draw[very thick] (0,0) circle(2.5cm);
      \draw[very thick, blue] (-0.25,-1.6) circle(0.5cm);
      \draw[very thick, red] (-1.1,1.2) circle(0.5cm);
      \draw[very thick, green] (1.4,0.8) circle(0.5cm);
    \end{scope}

    \begin{scope}[scale=0.3,xshift=4.5cm,yshift=1.2cm]
      \draw[very thick] (0,0) circle(2.5cm);
      \draw[very thick, blue] (-0.25,-1.6) circle(0.5cm);
      \draw[very thick, red] (-1.1,1.2) circle(0.5cm);
      \draw[very thick, green] (1.4,0.8) circle(0.5cm);
    \end{scope}
    
    \begin{scope}[scale=0.3,xshift=-0.5cm,yshift=-4cm]
      \draw[very thick] (0,0) circle(2.5cm);
      \draw[very thick, blue] (-0.25,-1.6) circle(0.5cm);
      \draw[very thick, red] (-1.1,1.2) circle(0.5cm);
      \draw[very thick, green] (1.4,0.8) circle(0.5cm);
    \end{scope}
  \end{tikzpicture}
  \caption{A Brunnian complex.}
  \label{fig:brunnian_complex}
\end{figure}

Furthermore, one might also require that the simplices have the
Brunnian property: remove one vertex, there are no other bonds.  For
example if we have $9$ vertices a Brunnian complex of type $(3,3)$
would look like the figure shown in Figure \ref{fig:brunnian_complex}.

$\HH$-complexes may be useful in describing situations where
non-pairwise interactions occur.

Finally, we suggest that hyperstructures may be useful in the analysis
of data --- extending hierarchical clustering to a hyperstructure
setting; see \citeasnoun{HTD}.  Bonds may also depend on parameters
and it may be useful to know the persistence (or stability) of the
bonds under parameter variations, for example for classification
purposes.

\section{Higher spaces}
\label{sec:spaces}

What is a space?  This is an old and interesting question.  We will
here add some higher (order) perspectives.  Often spaces are given by
open sets, metrics, etc.  They all give rise to bindings of points:
open sets, ``binding'' its points, distance binding points, etc.

In many contexts (of genes, neurons, links, subsets and subspaces,
$\ldots$) it seems more natural to specify the binding properties of
space by giving a hyperstructure --- even in addition to an already
existing ``space structure''.  In order to emphasize the binding
aspects of space we suggest that a useful notion of space should be
given by a set $X$ and a hyperstructure $\HH$ on it.  Such a pair
$(X,\HH)$ we will call a \emph{higher space}.  It tells us how the
points or objects are bound together.

Clearly there may be many such hyperstructures.  They may all be
collected into a larger hyperstructure --- $\HH^{\text{Total}}$ ---
which in a sense parametrizes the others.  Ordinary topological spaces
will be of order $0$ with open sets as bonds.  Through the bonds one
may now study processes like fusion and fission in the space.

Our key idea is that ``spaces'' and ``hyperstructures'' are intimately
connected.

In neuroscience one studies ``space'' through various types of cells:
place-, grid-, border-, speed-cells,$\ldots$; see \citeasnoun{Moser}.
All this spatial information should be put into the framework of a
``higher space'' with for example firing fields as basic bonds.  As
pointed out, the \emph{binding} problem fits naturally in here,
similarly ``cognitive'' and ``evolutionary'' spaces defined by
suitable hyperstructures.  Higher cognition should be described by
higher spaces.

Finally, the binding picture of space includes also the basic forces
of nature: gravity, electromagnetism, etc.\ which of course are
examples of bonds.  We have seen that links represent multiple bonds,
some of them not being reducible to pair bonds.

This raises again the natural question: \emph{Are there any basic
  forces in nature not reducible to pair bonds?  Which are the (basic)
  forces in nature (or artificial universes) that only act through
  clusters of objects (or higher bonds)?}  Think of quarks, Efimov
states and Brunnian links; see \citeasnoun{NS}.

\section{Conclusion}
\label{sec:conclusion}

We have discussed the philosophy and ideas around hyperstructures
which we argue captures the essence of organization and higher
structure.  We suggest that hyperstructures represent the form and
framework in which we think and act, and represent the external world
in our brains.

\subsection*{Acknowledgements}
I would like to thank A.~Stacey and M.~Thaule for help with graphics
and technical issues.  All the figures appear with their kind
permission.

\subsection*{Related links of interest:}
\begin{itemize}
  \item \href{http://www.newscientist.com/article/mg20927942.300-make-way-for-mathematical-matter.html}{New
    Scientist}\footnote{\href{http://www.newscientist.com/article/mg20927942.300-make-way-for-mathematical-matter.html}{\texttt{http://www.newscientist.com/article/mg20927942.300-make-way-for-
          mathematical-matter.htm}}}
  \item
    \href{http://www.technologyreview.com/view/422055/topologist-predicts-new-form-of-matter/}{MIT
      Tech
      Review}\footnote{\href{http://www.technologyreview.com/view/422055/topologist-predicts-new-form-of-matter/}{\texttt{http://www.technologyreview.com/view/422055/topologist-predicts-
          new-form-of-matter/}}}
\end{itemize}


\end{document}